\documentclass[amscd,amssymb,verbatim,12pt]{amsart}
\usepackage{floatflt}
\usepackage{hhline}
\usepackage{amssymb}
\usepackage{graphicx}
\usepackage[figuresright]{rotating}

\newfont{\myBbb}{msbm10 scaled 1200}
\newcommand{\N}{{\mbox{\myBbb N}}}
\newcommand{\Z}{{\mbox{\myBbb Z}}}

\begin{document}

\title{Exact properties of Frobenius numbers and fraction of the
symmetric semigroups in the weak limit for n=3}

\author{V. Shchur}
\address{Faculty of Mechanics and Mathematics,
Lomonosov Moscow State University,GSP-2, Leninskie Gory, Moscow,
119992, Russian Federation} \email{vladimir@chg.ru}


\date{\today}


\begin{abstract}
We generalize and prove a hypothesis by V. Arnold on the parity of
Frobenius number. For the case of symmetric semigroups with three
generators of Frobenius numbers we found an exact formula, which in
a sense is the sum of two Sylvester's formulaes. We prove that the
fraction of the symmetric semigroups is vanishing in the weak limit.
\end{abstract}
\maketitle

\section{Definitions}

Take $n$ mutually prime numbers $a_1, a_2, ..., a_n: (a_1, a _2,
..., a_n) = 1$. Consider semigroup $S = \left\{s=x_1a_1+...+x_na_n\|
x_i \in \Z_+\right\}$.\footnote{In the paper, $\N$ is the set of all
positive integers, and $\Z_+$ is the set of all non-negative
integers.} It follows that starting from some number $F(a_1, a_2,
..., a_n) \in S$ all integers are in the set $S$:

{\bf Definition 1}. $F(a_1, a_2, ..., a_n) = \min(s\in S\| \forall
k\in\Z_+, k\geq s: k\in S)$.

Number $F(a_1, a_2, ..., a_n)$ is called a Frobenius
number~\footnote{We have to note that often the Frobenius number
defined as a number less by unity (compare with the function $C$
defined below). In the paper we will keep given notation which is
also used in papers by V.Arnold.} and it's properties are the
subject of the paper. For $n=2$ the exact formula is known, the
Sylvester formula~\cite{Sylvester}, for $F$: $F(a_1, a_2) = (a_1 -
1)(a_2 - 1)$. It is known~\cite{Curtis} that for an arbitrary
relatively prime set $(a_1,a_2,...,a_n)$, $n\ge 3$ the Frobenius
number cannot be expressed in terms of a finite set of polynomials.

There are two useful functions. Function $C(a_1, a_2, ..., a_n) =
F(a_1, a_2, .. ., a_n)-1$ is the maximal integer which is not in the
set $S$. Function $G(a_1, a _2, ..., a_n) = F(a_1, a_2, ..., a_n) -
1 + a_1 + a_2 + ... + a_n$ is an analog of the function $C$ for the
semigroup $S_1 = \left\{s = x_1a_1+...+x_na_n\| x_i \in \N\right\}$.

Note, that condition $s < F$ and $s \in S$ leads to the result that
$C-s \notin S$ (otherwise, $C \in S$). The opposite statement is not
true in the general case.

{\bf Definition 2}. Semigroup $S$ is called symmetric if $\forall s
< F, s \notin S: C-s \in S$.

Each semigroup generated by two elements is the symmetric one.

\section{Exact results}

Analyzing properties of Frobenius numbers and associated semigroups,
V. Arnold noticed that if two of three generating elements are even
then $F$ is even as well~\cite{seminar}. The following Theorem is
the generalization of the experimental observation and gives the
proof. Note, that renumbering of the generating elements does not
change the structure of the set $S$ and, in particular, does not
change value of $F$. For the sake of simplicity we shall use this
property in following without explicitly mentioning it.

{\bf Theorem 1}. Let us $(a_1, a_2, ..., a_{n-1}, a_n - 1) = d$.
Then $F$ is divisible by $d$.

{\sf Proof}. By the definition of $F$: $F-1+a_n \in S (a_n \ge 1)$.
Hence, there exist nonnegative integers $x_1, ..., x_n$ such that
$F-1+a_n = x_1a_1 + \ldots + x_na_n$. Therefore, $F-1 = x_1a_1 +
\ldots + x_{n-1}a_{n-1} + (x_n-1)a_n$. If $x_n \ge 1$, then $F-1 \in
S$, what contradicts to the definition of $F$. Hence, $x_n = 0$ и $F
= x_1a_1 + \ldots + x_{n-1}a_{n-1} - (a_n - 1)$. Thus, $F$ is
divisible by $d$. $\blacksquare$

 Let $d_{ij} = (a_i, a_j), i \neq j$.

{\bf Lemma 1}. Johnson's lemma~\cite{Johnson}. $G(a_1, a_2, a_3) =
d_{12} G\left(\frac{a_1}{d_{12}}, \frac{a_2}{d_ {12}}, a_3\right)$.

{\bf Lemma 2}. Brauer-Shockley lemma~\cite{Brauer}. Let $d =
(a_1,a_2, ... , a_{n-1})$. Then, $C(a_1,a_2,..., a_n) = dC(a_1/d,
... , a_{n-1}/d, a_n)+ (d - 1)a_n$.

Divisibility of $F$ by $d$ follows also~\footnote{We acknowledge to
referee for this remark.} from the Lemma 2. Indeed, $F(a_1, . . . ,
a_n)=C(a_1, . . . , a_n) +1 = dC(a_1/d, . . . , a_{n-1}/d, a_n)+ (d
- 1)a_n + 1 = d(C(a_1/d, . . . , a_{n-1}/d, a_n)+ a_n)-(a_n-1)$.

{\bf Lemma 3}. Semigroup $S(a_1, a_2, a_3)$ is symmetric if and only
if for some $i$: $a_i \in S\left(\frac{a_j}{d_{jk}},
\frac{a_k}{d_{jk}}\right)$, where $i\neq j\neq k$.

In the case of symmetric semigroups we can prove the following exact
formula for the Frobenius numbers with $n=3$. Denote
$b_l=\frac{a_l}{d_{12}}, l=1, 2$. It is always possible to
renumerate generating elements in such a way that under the
conditions of Lemma 2 it takes place that $i=3$.

{\bf Theorem 2}. Let $S(a_1, a_2, a_3)$ is a symmetric semigroup,
$a_3 \in S(b_1, b_2)$. Then $F(a_1, a_2, a_3) =
d_{12}F\left(\frac{a_1}{d_{12}}, \frac{a_2}{d_{12}}\right) + F(d_
{12}, a_3)$.

{\sf Proof}. $F(a_1, a_2, a_3) = G(a_1, a_2, a_3) - a_1 - a_2 - a_3
+ 1 = $/by Lemma 1/ $= d_{12}G\left(\frac{a_1}{d_{12}},
\frac{a_2}{d_{12}}, a_3\right) - a_1 - a_2 - a_3 +1 = d_{ 12}(F(b_1,
b_2, a_3) + b_1 + b_2 + a_3 - 1) - d_{12}b_1 - d_{12}b_2 - a_3 +1 =
d_{ 12}F(b_1, b_2) + d_{12}a_3 - a_3 - d_{12} + 1 = $/by Sylvester's
formula/$ = d_{12}F (b_1, b_2) + F(d_{12}, a_3)$. Note, that
$(d_{12}, a_3) = 1$ due to the $(a_1, a_2 , a_3) = 1$.
$\blacksquare$

{\bf Notice}. If semigroup $S(a_1, a_2, a_3)$ is not symmetric then
$F(a_ 1, a_2, a_3) < d_{12}F\left(\frac{a_1}{d_{12}},
\frac{a_2}{d_{12}}) + F(d_{12}, a_3\right)$ , because if $a_3 \notin
S(b_1, b_2)$ then $F(b_1, b_2, a_3) < F(b_1, b_2)$.

At the time of the paper revision the author found in the literature
theorem in Chinese language~\cite{Chinese} similar to our Theorem 2.

\section{Asymptotic properties}

Theorem 2 shows that, informally speaking, if $a_3$ is sufficiently
large ($a_3 \ge  F(b_1, b_2)$), than $F(a_1, a_2, a_3)$ behaves as a
Frobenius number with two generating elements (compare with the
Sylvester formula!). Therefore, it is interesting to consider the
case under which all numbers $a_1, a_2, a_3$ are of the same order
in some sense. V. Arnold~\cite{Arnold} proposed to investigate an
asymptotic behavior in the following way. Let us fix vector ${\bf a}
= (a_1, ..., a_n)$ and for some function $f({\bf a})$ consider the
following average:

 $$V(f, {\bf a} r, N) = \sum_{-r \leq r_i \leq r, i = 1, ..., n} f(N {\bf a}+{\bf
 r}) / (2r)^n,$$ where ${\bf r} = (r_1, ..., r_n)\in \Z^n$ and $r, N \in \N.$
 Let us choose set $r(N)$ such that
\begin{equation}r(N) \to \infty, \frac{r(N)}{N}\to 0 \;\; {\rm if}\;\; N \to
\infty. \label{condition}
\end{equation}
Let us consider the limit:
 $$L(f, {\bf a}) = \lim_{r\to\infty, N\to\infty, r/N\to0} (V).$$

Let $Sym({\bf a}) = \left\{\begin{array}{llr} 1 & {\rm if} & S(a_1,
a_2, a_3)\; - \; {\rm symmetric} \\ 0 & &  - \; {\rm
otherwise}\end{array}\right.$. Fel~\cite{Fel} gives theorem which in
our terminology sounds as follows.

{\bf Theorem 3}. $L(Sym, {\bf a}) = 0$.

Unfortunately, proof given in~\cite{Fel} contains some mistakes. For
example, author says that set of triples $(a_1N+r_1, a_2N+r_2,
a_3N+r_3), |a_i| \le r(N), i=1,2,3$, which constitute symmetric
semigroups and such that $D=(a_1N+r_1, a_2N+r_2) \sim N$, is empty
while $N \to \infty$. Let us give counter example. Let $a_1 = 4, r_1
= 2, a_2 = 6, r_2 = 3$. Then $D=2N+1\sim N$ and
$F(\frac{a_1N+r_1}{D}$, $\frac{a_2N+r_2}{D}) = 2$, i.e. for any set
$r(N)$, fulfilling~(\ref{condition}), and any $a_3$ for large enough
$N$ semigroup $S(a_1N+r_1, a_2N+r_2, a_3N+r_3)$ (if exists) is
symmetric. Thus, considered set are not always empty. We present
here some other proof of the Theorem 3. The main idea is not to
prove an emptiness but rather prove the fact that the fraction of
symmetric semigroups is small in the considered limit.

We say that the fraction of the symmetric semigroups is vanishing in
the weak asymptotics. Fel showed~\cite{Fel} that in the limit of
interest the fraction of the mutually simple numbers (those sets on
which Frobenius numbers defined) is equal $\frac{1}{\zeta(n)}$. In
other words for

$I({\bf a}) = \left\{
\begin{array}{llr}
1 & {\rm if} & (a_1, ..., a_n) = 1\\
 0 & & - \; {\rm otherwise}\end{array}\right.$
we have $L(I, {\bf a}, 1) = \frac{1}{\zeta(n)}$.

{\bf Proof}. Let $A_i = N a_i + r_i$. Let us figure out the number
of semigroups in the cubic $r$-vicinity of the point $N{\bf a}$
those for which $A_3 \in S_0 = S\left(\frac{A_1}{D},
\frac{A_2}{D}\right)$, where $D = (A_1, A_2)$. Note, that for
sufficiently large $D$ Frobenius number
$F\left(\frac{A_1}{D},\frac{A_2}{D}\right)$ is not larger then
$A_3$, therefore for sufficiently large $D$ any $A_3$ lies in
semigroup
 $S_0$. The number of semigroups such that $D \geq D_0$ (value of $D_0$ will
 be chosen later on) may be estimated as:
$$
U_1 = \sum_{D = D_0}^{cN}\left(\frac{2r}{D}\right)^2 \leq
4r^2\int\limits_{D_0-1}^ {cN-1} \frac{1}{x^2}dx =
4r^2\left(\frac{1}{D_0-1}-\frac{1}{cN}\right),
$$
where constant $c$ may be chosen as $c=\min(a_1, a_2)+1$ because if
$N\to\infty , r\to\infty, r/N\to0$ then $D\leq \min(A_1, A_2)\leq
cN$.

Now consider the case $D < D_0$. Let us fix some $\varepsilon>0$.
Then $\exists N_0$ such that $\forall N>N_0$ we have
$r(N)<\varepsilon N$, i.e. $\forall N>N_0$ $\forall r_3$:
$|r_3|<r(N)$ we have that $a_3N+r_3<(a_3+\varepsilon)N$. Denote
$c_1=a_3+\varepsilon$. Let us estimate the number of elements of the
semigroup $S\left(\frac{A_1}{D},\frac{A_2}{D}\right)$ in the
interval $[0;c_1N]$: $s=x_1a_1+x_2+a_2$. (See Fig.~1) Without lack
of generality let that $a_1<a_2$. Then for sufficiently large $N$
$\frac{A_1}{D}<\frac{A_2}{D}$. Denote $T =
\left[\frac{c_1ND}{A_1}\right]$. Notice that $x_1+x_2\le T$. With
fixed number $x_2$ number $x_1$ may take not more than  $T-x_2$
values. Therefore,
$$\#\left(S\left(\frac{A_1}{D},\frac{A_2}{D}\right)\bigcap\left[0;c_1N\right]\right)\le
\sum^{T}_{x_2=0}\left(T-x_2\right)=
\sum^{T}_{k=0}k=\frac{T\left(T+1\right)}{2}.$$ Then with fixed $A_1$
and $A_2$ and any $A_3=a_3N+r_3$, $|r_3|\le r(N)$ the number of
symmetric semigroups are not larger than
$\left(S\left(\frac{A_1}{D},\frac{A_2}{D}\right)\bigcap\left[0;c_1N\right]\right).$
So the number of symmetric semigroups such that $D<D_0$, may be
estimated as
$$
U_2\le\sum_{D=1}^{D_0}\left(\frac{2r}{D}\right)^2\left(C_1^2\left(\frac{N}{A_1}\right)^2D^2+O(D)\right)\approx
4r^2C_1^2\left(\frac{N}{a_1N+r_1}\right)^2D_0.
$$

Thus, the fraction of symmetric semigroups is estimated as
\begin{equation}
W = 3\frac{U_1 + U_2}{8r^3} \leq 3\left(\frac{1}{D_0-1}-\frac{1}{cN}
+ \frac{D_0}{2r}C_1^2\left(\frac{N}{a_1N+r_1}\right)^2 \right).
\label{exp-W}
\end{equation}

Let $D_0 = \ln r$. Then, with $N \to \infty, r \to \infty, r/N \to
0$ the right hand side of inequality (\ref{exp-W}) vanishes, i.e. $W
\to 0$. $\blacksquare$

\medskip

{\bf Acknowledgment}. Author acknowledge scientific advises by V.I.
Arnold and S.M. Gusein-Zade and help from A. Ustinov. This work
supported by the Russian Federation Grant NSh-709.2008.1.

\newpage

\begin{figure}
\centering
\includegraphics[angle=0,width=\columnwidth]{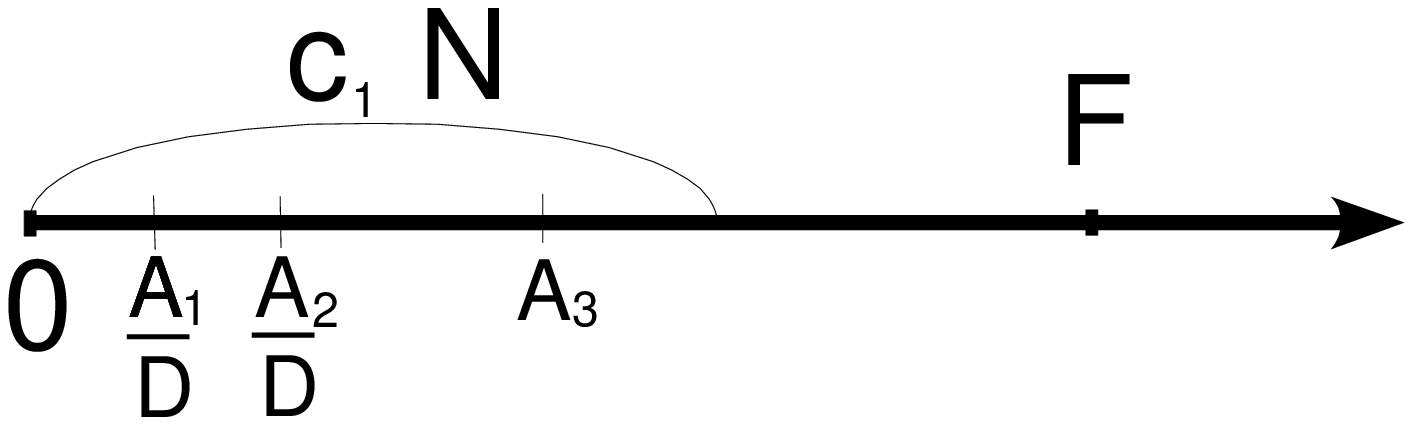} \vspace{-13cm} \caption{}\label{fig1}
\end{figure}

\end{document}